\documentstyle[12pt]{article}
\begin{document}
\newpage
\pagestyle{empty}
\setcounter{page}{0}
%
%%% ** start of amsfont definitions **
%\newfont{\twelvemsb}{msbm10 scaled\magstep1}
%\newfont{\eightmsb}{msbm8} \newfont{\sixmsb}{msbm6} \newfam\msbfam
%\textfont\msbfam=\twelvemsb \scriptfont\msbfam=\eightmsb
%\scriptscriptfont\msbfam=\sixmsb 
%\def\Bbb{\ifmmode\let\next\Bbb@\else \def\next{\errmessage{Use
%      \string\Bbb\space only in math mode}}\fi\next}
%\def\Bbb@#1{{\Bbb@@{#1}}} \def\Bbb@@#1{\fam\msbfam#1}
%4\newfont{\twelvegoth}{eufm10 scaled\magstep1}
%\newfont{\tengoth}{eufm10} \newfont{\eightgoth}{eufm8}
%\newfont{\sixgoth}{eufm6} \newfam\gothfam
%\textfont\gothfam=\twelvegoth \scriptfont\gothfam=\eightgoth
%\scriptscriptfont\gothfam=\sixgoth \def\frak{\frak@}
%\def\frak@#1{{\fam\gothfam{{#1}}}} \def\frak@@#1{\fam\gothfam#1}}
%%% ** end of amsfont definitions **
%
%
%
%\def\Bbb{\bf}
\def\CC{{\Bbb C}}
\def\NN{{\Bbb N}}
\def\QQ{{\Bbb Q}}
\def\RR{{\Bbb R}}
\def\ZZ{{\Bbb Z}}
\def\cA{{\cal A}}          \def\cB{{\cal B}}          \def\cC{{\cal C}}
\def\cD{{\cal D}}          \def\cE{{\cal E}}          \def\cF{{\cal F}}
\def\cG{{\cal G}}          \def\cH{{\cal H}}          \def\cI{{\cal I}}
\def\cJ{{\cal J}}          \def\cK{{\cal K}}          \def\cL{{\cal L}} 
\def\cM{{\cal M}}          \def\cN{{\cal N}}          \def\cO{{\cal O}}
\def\cP{{\cal P}}          \def\cQ{{\cal Q}}          \def\cR{{\cal R}} 
\def\cS{{\cal S}}          \def\cT{{\cal T}}          \def\cU{{\cal U}}
\def\cV{{\cal V}}          \def\cW{{\cal W}}          \def\cX{{\cal X}}
\def\cY{{\cal Y}}          \def\cZ{{\cal Z}}
\def\qed{\hfill \rule{5pt}{5pt}}
\newtheorem{lemma}{Lemma}
\newtheorem{prop}{Proposition}
\newtheorem{theo}{Theorem}
\newenvironment{result}{\vspace{.2cm} \em}{\vspace{.2cm}}

\rightline{CPTH-9611476}  
\rightline{q-alg/9611014} 
\rightline{November 96} 

\vfill
\vfill

\begin{center}

{\LARGE {

\bf {\sf 
Induced Hopf stucture and irreducible representations of an elliptic 
${\cal U}_{q,p}(sl(2))$ via a nonlinear map
}}}\\[2cm]

\smallskip 

{\large A. Chakrabarti\footnote{chakra@orphee.polytechnique.fr}}

\smallskip 

\smallskip 

\smallskip

{\em  \footnote{Laboratoire Propre du CNRS UPR A.0014}Centre de Physique 
Th\'eorique, Ecole Polytechnique, \\
91128 Palaiseau Cedex, France.}

\end{center}

\vfill

\begin{abstract}
Shiraishi's two parameter generalization of ${\cal U}_{q}(sl(2))$ to 
${\cal U}_{q,p}(sl(2))$ involving an elliptic function is considered. 
The generators are mapped non-linearly on those of ${\cal U}_{q}(sl(2))$. 
This gives directly the irreducible representations and an induced Hopf 
structure. This is one particular example of the scope of a class of 
non-linear maps introduced by us recently.  
\end{abstract}

\vfill
\vfill

\newpage
\pagestyle{plain}

Recently Shiraishi [1] has made what he describes as " an attempt at obtaining
an elliptic $sl(2)$ algebra...". He generalizes the ${\cal U}_{q}(sl(2))$ 
algebra as follows
\begin{eqnarray}
&& q^{2J_{0}}{\hat J}_{\pm}q^{-2J_{0}}= q^{\pm 2}{\hat J}_{\pm}  \\
&& [{\hat J}_{+},{\hat J}_{-}]=\sum_{n\in Z}(-1)^{n}q^{2J_{0}(2n+1)}
p^{(n+1/2)^{2}}
\end{eqnarray}
where for his $(t,e,f)$ we have written ($q^{2J_{0}},{\hat J}_{+},
{\hat J}_{-}$) respectively. The generalization involves a theta function on 
the rhs of (2) insead of $[2J_{0}]$ (i.e.
$(q^{2J_{0}}-q^{-2J{0}})(q-q^{-1})^{-1})$ giving standard ${\cal U}_q(sl(2))$. This is
proposed in the context of the  "elliptic algebra" ${\cal A}_{q,p}({\hat sl(2)})$
of [2]. The connection between the  two formalismes ([1] and [2]) is not clear.
But the algebra given by (1) and (2) (which will be denoted by ${\cal U}_{p,q}(sl(2))$
has some interesting  properties. A Heisenberg-Clifford realization is presented
in [1].

We point out in this note that non-linear mappings of [3] can be adapted to
provide one between ${\cal U}_{q}(sl(2)$ and ${\cal U}_{p,q}(sl(2))$. Such a 
map immediately gives the irreducible representations of ${\cal U}_{p,q}
(sl(2))$ and provides an induced Hopf stucture(absent in [1]).

Let us recapitulate the formalisme of Sec.3 of [3] in a form well-adapted
to the present case. 
The generators of ${\cal U}_{q}(sl(2))$ satisfy
\begin{eqnarray}
&& q^{2J_{0}}J_{\pm}q^{-2J_{0}} =q^{\pm 2}J_{\pm}  \\
&&[J_{+},J_{-}] =[2J_{0}]
\end{eqnarray}
The Casimir operator is
\begin{eqnarray}
&& C= J_{-}J_{+}+[J_{0}][J_{0}+1]
\end{eqnarray}
Define (introducing a function $\phi$ with suitable properties to be 
specified below),
\begin{eqnarray}
&& {\hat J}_{+}=J_{+}\biggl( { \phi(C) - \phi([J_{0}][J_{0}+1]) \over
C -[J_{0}][J_{0}+1]}\biggr)^{1+ \eta \over 2} \\
&& {\hat J_{-}}= \biggl({\phi(C)-\phi([J_{0}][J_{0}+1])\over
C-[J_{0}][J_{0}+1]}\biggr)^{1-\eta\over 2} J_{-}
\end{eqnarray}
where $\eta=0,1,-1.$ (The choice $\eta=0$ will be the standard one. But in 
certain contexts it might be prefable to avoid squareroots by setting $\eta=
1$ or $-1$. In the latter cases the familiar matrix elements of$J_{\pm}$
should also be consistently written as $([j][j+1]-[m][m\pm 1])^{1\pm\eta\over
2}$.)
 
It follows easily (for any $\eta$ and using (5)) that
\begin{eqnarray}
&& [{\hat J}_{+},{\hat J}_{-}]=\phi([J_{0}][J_{0}+1])-\phi([J_{0}][J_{0}-1])
\end{eqnarray}
Hence for any suitably chosen function $\chi(J_{0})$ to ensure
\begin{eqnarray}
&&[{\hat J}_{+},{\hat J}_{-}]=\chi(J_{0})
\end{eqnarray}
one must have
\begin{eqnarray}
&&\phi([J_{0}][J_{0}+1])- \phi([J_{0}][J_{0}-1]=\chi(J_{0})
\end{eqnarray}
This is a general result. In [3] we gave the example that for
\begin{eqnarray}
&&\chi(J_{0})= [2J_{0}](1+\beta[2] [J_{0}]^2) \\
&&\phi(x)=x+ \beta x^2
\end{eqnarray}
For the general case, it is convenient to define
\begin{eqnarray}
&&\phi([J_{0}][J_{0}+1])= \psi(J_{0})= \sum_{n\in Z} a_nq^{2nJ_{0}}
\end{eqnarray}
limiting our considerations to such a series in  $t^{\pm1}=q^{\pm 2J_{0}}$.
It is not essential to reconvert $\psi(J_{0})$ into $\phi([J_{0}][J_{0}+1])$.

Similarly one may define (in (6) and (7))
\begin{eqnarray}
&& \phi(C)= \psi(J_{op})
\end{eqnarray}
where
\begin{eqnarray}
&& C=[J_{op}][J_{op}+1]
\end{eqnarray}
giving $q^{2J_{op}}$ as the solution of a quadratic equation.

Let in (9) and (10)
\begin{eqnarray}
&& \chi(J_{0}) = \sum_{k\geq1} ( b_{k}q^{2kJ_0} +b_{-k}q^{-2kJ_0})
\end{eqnarray}

From (10), (13) and (16)
\begin{eqnarray}
&&
a_{k}={q^kb_k\over(q^k-q^{-k})},\qquad a_{-k}=-{q^{-k}b_{-k}\over(q^k-q^{-k})}\\
&& (k=1,2,\cdots) \end{eqnarray}
while $a_0$ in (13) is arbitrary (and indeed cancels out in the difference
$(\psi(J_{op})-\psi(J_0))$.
Hence 
\begin{eqnarray}
&&\psi(J_0)= a_0+ \sum_{k\geq1}{(b_kq^{k(2J_0+1)}-b_{-k}q^{-k(2J_0+1)})\over
(q^k-q^{-k})}
\end{eqnarray}
Choosing, for example, $c_0$ being a finite constant,
\begin{eqnarray}
&&a_0=c_0 -\sum_{k\geq1}{(b_k-b_{-k})\over(q^k-q^{-k})}
\end{eqnarray}
$\psi(J{_0})$ can be, separately, given a finite limit as $ q\rightarrow 1$.
(But this is not strictly necessary as $a_0$ cancels in the numerator.)

Writing (12) in the form (13) (11 in the form(16)) one has
\begin{eqnarray}
&&a_{\pm 1}={q^{\pm 1}\over(q-q^{-1})^2}(1-{2\beta\over(q-q^{-1})^2})
= \pm{q^{\pm 1}\over(q-q^{-1})}b_{\pm 1}\\
&&a_{\pm 2}={q^{\pm 2}\over(q-q^{-1})^4}{\beta\over(q+q^{- 1})}
 =\pm{q^{\pm 2}\over(q^2-q^{-2})}b_{\pm 2}
\end{eqnarray}

Corresponding to (2) one has
\begin{eqnarray}
&&\chi(J_0)= \sum_{n\in Z}(-1)^nq^{2j_0(2n+1)}p^{(n+{1\over 2})^2} \nonumber\\
&&\phantom{\chi(J_0)}=\sum_{k\geq 1}(b_kq^{2kJ_0}+b_{-k}q^{-2kJ_0})
\end{eqnarray}
with for even and odd $k$ respectively
\begin{eqnarray}
&&b_{\pm k}=0,\qquad b_{\pm k}=(-1)^{k\pm 1\over 2}p^{({k\over 2})^2}.
\end{eqnarray}
Now (17) gives the corresponding $a_{\pm k}$'s.

 Considering generic $q$, our map gives, quite generally, the $(2j+1)$ 
dimentional irreducible representations (for (half) integer $j$)
\begin{eqnarray}
&& q^{\pm 2J_0}\mid j,m\rangle= q^{\pm 2m}\mid j,m\rangle\qquad 
(m=j,j-1,\cdots,-j)\\
 &&{\hat J}_{\pm} \mid j,m\rangle =(\psi(j)-\psi(m))^{1\pm
\eta\over 2}\mid j,m\pm 1\rangle\\
 &&\phantom{{\hat J}_{\pm} \mid j,m\rangle}=
\biggl(\phi([j][j+1])-\phi([m][m+1]\biggr)^{{1\pm \eta}\over 2}\mid j,m\pm 1\rangle\\
&& \qquad\qquad\qquad \qquad\qquad           (\eta= 0,1,-1)
\end{eqnarray}

It is convenient to express the Casimir now as
\begin{eqnarray}
&& {\hat C}= {\hat J_-}{\hat J_+} + \phi([J_0][J_0+1])
\end{eqnarray}
so that
\begin{eqnarray}
&&{\hat C}\mid j,m\rangle= \phi([j][j+1])\mid j,m\rangle\\
&&\phantom{ {\hat C}\mid j,m\rangle } =\psi(j)\mid j,m\rangle
\end{eqnarray}

For ${\cal U}_q(sl(2))$ one has the coproducts
\begin{eqnarray}
&&\Delta J_0=J_0 \otimes 1 +1 \otimes J_0\\
&&\Delta J_{\pm} =J_{\pm}\otimes q^{J_0} +q^{-J_0} \otimes J_\pm\\
&&\Delta C = (\Delta J_-)(\Delta J_+) +[\Delta J_0][\Delta J_0+1]
\end{eqnarray}
They induce the coproducts
\begin{eqnarray}
&& \Delta{\hat J}_{+}=\Delta J_{+}\biggl( { \phi(\Delta C) -
\phi([\Delta J_{0}][\Delta J_{0}+1]) \over\Delta C -[\Delta
J_{0}][\Delta J_{0}+1]}\biggr)^{1+ \eta \over 2} \\
 && \Delta{\hat J_{-}}=
\biggl({\phi(\Delta C)-\phi([\Delta J_{0}][\Delta J_{0}+1])\over
(\Delta C)-[\Delta J_{0}][\Delta J_{0}+1]}\biggr)^{1-\eta\over 2}\Delta J_{-}
\end{eqnarray}

The counits and the antipodes are analogously induced. The results from (25)
to (36) are general. For each particular case one has to implement the
appropriate $\psi$ or $\phi$ as discussed above.

In (25) to (27) we have considered $(2j+1)$ dimentional irreps for generic $q$.
But our formalism can be implemented also for $q$ a root of unity. The 
parameters $j$ and $m$ in (26) can have "fractional parts" [4] and one can 
have, for example,
periodic and semiperiodic representations of entirely different dimentions. 
But we will not study these aspects here. The infinite series form in (2) is
not appropriate for spaces of periodic representations. Here our purpose has
been to draw attention to the fact that irreps and an induced Hopf structure
corresponding to the type given by (2) are quite simply furnished by our
nonlinear map.

A complementary class of nonlinear map has been introduced [5] relating
${\cal U}(sl(2))$ and ${\cal U}_h(sl(2))$. It is complementary in the sense 
that here (for ${\cal U}(sl(2)$ or ${\cal U}_q(sl(2))$ as starting point) 
the nonlinearity is implemented via the diagonalizable generator $J_0$ (or 
$q^{\pm J_0}$) whereas in [5] the corresponding role is played by the 
nondiagonalizable generator $J_+$. These two types of mappings can be 
combined. Such a study will be presented elsewhere.

\end{document}